\documentclass[conference]{IEEEtran}

\usepackage{amsfonts}
\usepackage{cite,url,subfigure,epsfig,graphicx}
\usepackage{amssymb,amsmath}

\usepackage{balance}

\ifCLASSINFOpdf

\else

\fi

\IEEEoverridecommandlockouts

\begin{document}
\title{
Adaptive Robust Energy Management Strategy for Campus-Based Commercial Buildings Considering Comprehensive Comfort Levels
}

\author{\IEEEauthorblockN{Zheming Liang\IEEEauthorrefmark{2}\IEEEauthorrefmark{4},
Desong Bian\IEEEauthorrefmark{2},
Dawei Su\IEEEauthorrefmark{3},
Ruisheng Diao\IEEEauthorrefmark{2},
Di Shi\IEEEauthorrefmark{2}\IEEEauthorrefmark{1},
Zhiwei Wang\IEEEauthorrefmark{2} and
Wencong Su\IEEEauthorrefmark{4}
}
\IEEEauthorblockA{\IEEEauthorrefmark{2}GEIRI North America, San Jose, CA, USA}
\IEEEauthorblockA{\IEEEauthorrefmark{3}State Grid Jiangsu Electric Power Company, Nanjing, China}
\IEEEauthorblockA{\IEEEauthorrefmark{4}Department of Electrical and Computer Engineering, University of Michigan-Dearborn, Dearborn, MI, USA}
\thanks{This work is funded by SGCC Hybrid Energy Storage Management Platform for Integrated Energy System.

Di Shi is the corresponding author. (Email: di.shi@geirina.net)
}
}

\maketitle

\begin{abstract}
Neglecting consumers' comfort always leads to failure or slow-response to demand response request. In this paper, we propose several comprehensive comfort level models for various appliances in campus-based commercial buildings (CBs). The objective of the proposed system is to minimize O\&M costs of campus-based CBs and maximize various comfort levels simultaneously under the worst-case scenarios. Adaptive robust optimization (ARO) is leveraged to handle various uncertainties within the proposed system: (i) demand response signals sending from the distribution system operator (DSO); (ii) arrival state-of-charge (SoC) conditions of plug-in electric vehicles (PEVs); (iii) power outputs of renewable energy sources (RESs); and (iv) load demand of other appliances. Benders decomposition, such as column-and-constraint generation (C\&CG) algorithm, is used to solve the reformulated NP-hard min-max problem. Extensive simulation results demonstrate the effectiveness of the proposed optimal energy management strategy for campus-based CBs in both minimizing O\&M costs and maximizing comprehensive comfort levels.
\end{abstract}

\IEEEpeerreviewmaketitle

\section{Introduction}\label{section1}
With the development of economy all over the world, new appliances such as plug-in electric vehicles (PEVs), heating, ventilation, air-conditioning (HVAC) systems, and roof-top solar panels are vastly deployed in distribution systems which increase the total power demand significantly, especially for peak hours. However, capital costs of installing new generators and enhancing distribution lines and substations are too expensive compared with reducing and shifting peak loads. Therefore, demand response (DR) is introduced as a new method to shift the peak loads to off-peak hours or curtail non-critical loads directly to avoid contingency issues on substations and transmission/distribution lines. Based on objectives, prior works can be categorized into three different types, namely, (i) minimizing operation and maintenance (O\&M) costs of the entire system~\cite{LiAl17}, (ii) maximizing consumers' comfort levels (only focus on the HVAC system)~\cite{KoBa16}, and (iii) minimizing load curtailment costs~\cite{KhSh11}. Very few of them focus on joint objectives of both minimizing O\&M costs and maximizing consumers' comfort levels~\cite{LiBi19}.

More importantly, commercial buildings (CBs) consume more than $40\%$ of total energy supply to power systems~\cite{KlKw12}. Compared with residential households, campus-based CBs generally have larger peak demand and more total energy consumption. Therefore, CBs are better venues in performing DR to enhance the stability of the distribution system. However, as aforementioned, consumers' comfort levels are usually sacrificed when optimal DR results are achieved, which is typically caused by: (i) shutting down HVAC systems in summer which increases indoor temperatures and vice versa; (ii) shutting down electric water heaters (EWHs) which decreases hot water temperatures; and (iii) stop charging PEVs before their batteries are fully charged~\cite{BiPi15}. Therefore, including consumers'/occupants' comfort levels will boost the adoption rates of various DR programs.

Additionally, the amount of aforementioned new appliances installed in CBs bring a lot of uncertainties into the proposed system~\cite{LiGu17}. For instance, arrival state-of-charges (SoCs) of PEVs are unknown, making it difficult for central controllers of CBs to obtain a probability density function (PDF) of PEV related uncertainty. However, a range including upper and lower bounds for arrival SoCs of PEVs is easy to obtain. Similarly, for other uncertainties, such as: (i) power outputs of renewable energy sources (RESs), e.g., roof-top solar panels; (ii) elastic base load demand; and (iii) demand response signal from the distribution system operator (DSO), ranges can be acquired. Therefore, an adaptive robust optimization (ARO) approach can be leveraged to handle these uncertainties to ensure optimal energy management strategies under worst-case scenarios. Cardinality uncertainty sets are generated based on ranges of these uncertainties. The ARO approach can overcome the over-conservativeness issue through selecting proper budgets of uncertainty based on uncertainty sets. Besides, the ARO approach requires far less computational time compared with scenario-based stochastic optimization (SO) approach~\cite{LiCh17}.

In this paper, we propose several comprehensive comfort level models for various appliances in campus-based CBs. The objective of our proposed system is to minimize O\&M costs of campus-based CBs while maximize various comfort levels simultaneously under worst-case scenarios. ARO is leveraged to handle various uncertainties within the proposed system. Benders decomposition, such as column-and-constraint generation (C\&CG) algorithm, is used to solve the reformulated NP-hard min-max problem. Extensive simulation results demonstrate the effectiveness of our proposed optimal energy management strategy for campus-based CBs.

The following contributions are made in this paper:
\begin{itemize}
\item Models for comfort levels of appliances in the campus-based CBs are proposed. The objective is to minimize CBs' O\&M costs while maximize consumers' comfort levels simultaneously.

\item An ARO approach is introduced to handle various uncertainties. C\&CG algorithm is used to solve the reformulated NP-hard min-max problem.

\item Extensive simulation results show the effectiveness of the proposed framework and the proposed comfort models are proved to fit a novel energy management strategy.
\end{itemize}

The rest of this paper is organized as follows. Section~\ref{section2} presents the mathematical formulation of the optimal energy management problem. Section~\ref{section3} proposes an efficient solution method to solve the optimization problem. Section~\ref{section4} summaries the simulation results of our test system. Section~\ref{section5} draws our conclusions.

\section{Problem Formulation}\label{section2}
We consider the case that several CBs exist on one campus. Each CB consists of one HVAC system, one EWH, one electric storage system (ESS), one pack of roof-top solar panel, and a base power load (servers, lights, personal computers, projectors, etc.)~\cite{HaLi14}. Each CB has a parking lot with charging stations for PEVs. PEVs' penetration level in the system is defined as consumers with PEVs versus total populations in a CB. The campus is connected downstream of a distribution system through a substation. The operating horizon includes office hours of one operating day from 8 a.m. to 8 p.m. with 48 time slots.

\subsection{Constraints}

\subsubsection{Constraints of Appliances }
We have the following constraints to ensure the stability of the proposed system. Firstly, we have constraints for HVAC systems~\cite{LiBa15}:
\begin{equation}\label{con:HVAC1}
T_{i, t+1} = \beta_i T_{i, t} + \alpha_i U_{i, t}, \forall i, t,
\end{equation}
where $T_{i, t} = [T_{i, t}^{\text{in}}, T_{i, t}^{\text{iw}}, T_{i, t}^{\text{ow}}]^T$ includes indoor temperature, inner wall temperature, and outer wall temperature, respectively. Variable $U_{i, t} = [T_{t}^{\text{out}}, \Psi_t, \sigma_{i, t} \eta_i P_{i, t}^{\text{hvac}}]^T$ includes outdoor temperature, solar irradiance, cooling/heating indicator, and output cooling/heating of a HVAC system in CB $i$, respectively. Variables $\alpha_i$ and $\beta_i$ are environment coefficients of CB $i$~\cite{Th08}.
\begin{equation}\label{con:HVAC2}
T^d_{i} - \delta_{i} \leq T_{i, t}^{\text{in}} \leq T^d_{i} + \delta_{i}, 0 \leq P_{i, t}^{\text{hvac}} \leq \overline{P}_{i}^{\text{hvac}}, \forall i, t,
\end{equation}
where $T^d_{i}$ is the desired indoor temperature, $\delta_{i}$ is the maximum temperature deviation from the desired indoor temperature, $\overline{P}_{i}^{\text{hvac}}$ is the maximum power consumption of a HVAC system. Moreover, consumer's comfort level related to a HVAC system in the $i$-th CB can be defined as:
\vspace{-5pt}
\begin{equation}
J_{i, t} = \begin{cases}\label{con:comfort1}
0, & T_{i, t}^{\text{in}} \geq T^d_{i} + \delta_{i} \\
\frac{T^d_{i} + \delta_{i} -T_{i, t}^{\text{in}}}{\delta_{i}-\epsilon_{i}}, & T^d_{i} + \epsilon_{i} \leq T_{i, t}^{\text{in}} \leq T^d_{i} + \delta_{i} \\
1, & T^d_{i} - \epsilon_{i} \leq T_{i, t}^{\text{in}} \leq T^d_{i} + \epsilon_{i} \\
\frac{T_{i, t}^{\text{in}} - (T^d_{i} - \delta_{i})}{\delta_{i}-\epsilon_{i}}, & T^d_{i} - \delta_{i} \leq T_{i, t}^{\text{in}} \leq T^d_{i} - \epsilon_{i}\\
0, & T_{i, t}^{\text{in}} \leq T^d_{i} - \delta_{i}.
\end{cases}\vspace{-5pt}
\end{equation}
Comfortable indoor temperature zone is defined as $T^d_{i} \pm \epsilon_{i}$, where $\epsilon_{i}$ is the maximum indoor temperature deviation from desired temperature that can still ensure a comfortable temperature zone.

Additionally, we have the following charging dynamics for PEVs:
\vspace{-8pt}
\begin{equation}\label{con:pev_dynamics1}
SoC_{v, t} = SoC_{v, 0} + \sum_{\tau=1}^t \frac{p_{v, \tau}^{\text{ch}} \eta_{v}^{\text{ch}}}{\overline{E}_v}, \forall v, t, \vspace{-5pt}
\end{equation}
where $p_{v, t}^{\text{ch}}$ is the charging rate, $\eta_{v}^{\text{ch}}$ is the charging efficiency, $\overline{E}_v$ is the rated energy.
\begin{equation}\label{con:pev_dynamics2}
\underline{SoC}_v \leq SoC_{v, t} \leq \overline{SoC}_{v}, 0 \leq p_{v, t}^{\text{ch}} \leq \overline{P}_{v}^{\text{ch}}, \forall v, t,
\end{equation}
where upper bound $\overline{SoC}_{v}$ and lower bound $\underline{SoC}_v$ are imposed to enhance batteries' lifetimes. Variable $\overline{P}_v^{\text{ch}}$ denotes maximum charged energy over period $t$ for the $v$-th PEV. Furthermore, comfort level related to the $v$-th PEV can be defined as follows:
\vspace{-15pt}
\begin{equation}
J_{v, t} = \begin{cases}
1, & SoC_{v}^d \leq SoC_{v, t} \\
\frac{SoC_{v, t} - SoC_{v}^{\text{base}}}{SoC^d_{v} - SoC_{v}^{\text{base}}}, & SoC_{v}^{\text{base}} \leq SoC_{v, t} \leq SoC_{v}^d \\
0, & SoC_{v, t} \leq SoC_{v}^{\text{base}}. \label{con:comfort2}
\end{cases}\vspace{-5pt}
\end{equation}
where $J_{v, t}$ denotes comfort level of the $v$-th PEV owner. Variable $SoC_{v}^d$ is the desired SoC for the $v$-th PEV. Variable $SoC_{v}^{\text{base}}$ represents the base SoC required for the $v$-th PEV with round trip between its home and a CB.

Moreover, the requirements for EWHs are expressed as:
\vspace{-8pt}
\begin{equation}\label{con:EWH3}
T_{l, t} = T_{l, 0} + \sum_{\tau=1}^{t} \frac{\zeta_l p_{l, \tau} - H_{l, \tau}^{\text{de}}}{M_l C_{\text{water}}}, \forall l, t, \vspace{-5pt}
\end{equation}
where $p_{l, \tau}$ is power consumption of the $l$-th EWH. Variable $H_{l, \tau}^{\text{de}}$ is heat decrease of the $l$-th EWH, including heat loss that is transferred to its ambient, outflow of hot water and inflow of cold water. Parameter $M_l$ is the mass of water in tank $l$, and $C_{\text{water}}$ is specific heat capacity of water. Furthermore, water temperature in the $l$-th EWH needs to be regulated within certain range to maintain comfort level related to EWHs:
\begin{equation}\label{con:EWH4}
T_{l}^d - \delta_{l} \leq T_{l, t}, \underline{P}_l \leq p_{l, t} \leq \overline{P}_l, \forall l, t
\end{equation}
\vspace{-15pt}
\begin{equation}
J_{l, \tau} = \begin{cases}\label{con:comfort3}
1, & T_{l}^d \leq T_{l, t} \\
\frac{T_{l, t} - (T_{l}^d - \delta_{l})}{T^d_{l} - (T_{l}^d - \delta_{l})}, & T_{l}^d - \delta_{l} \leq T_{l, t} \leq T_{l}^d \\
0, & T_{l, t} \leq T_{l}^d - \delta_{l}.
\end{cases}\vspace{-5pt}
\end{equation}
where $T^d_{l}$ is the desired water temperature in the $l$-th EWH. Variable $\delta_{l}$ is the maximum allowed temperature deviation from the desired water temperature.

Similarly, we have the following dynamics for ESSs:
\vspace{-8pt}
\begin{equation}\label{con:es_dynamics1}
SoC_{k, t} = SoC_{k, 0} + \sum_{\tau=1}^t \frac{p_{k, \tau}^{\text{ch}} \eta_{k}^{\text{ch}}}{\overline{E}_k} - \sum_{\tau=1}^t \frac{p_{k, \tau}^{\text{dis}}/\eta_{k}^{\text{dis}}}{\overline{E}_k}, \forall k, t,
\end{equation}
where $p_{k, t}^{\text{ch}}$ and $p_{k, t}^{\text{dis}}$ are power charged into or discharged from the $k$-th ESS at time $t$. Variables $\eta_{k}^{\text{ch}}$ and $\eta_{k}^{\text{dis}}$ represent charging and discharging efficiencies of the $k$-th ESS, respectively. Each ESS has a finite capacity, therefore, the energy stored in it must has the following lower and upper bounds:
\begin{equation}\label{con:es_dynamics2}
\underline{SoC}_k \leq SoC_{k, t} \leq \overline{SoC}_{k}, \forall k, t,
\end{equation}
where $\overline{SoC}_{k}$ is the upper bound and $\underline{SoC}_k$ is the lower bound of the $k$-th ESS' SoC. Furthermore, ES units have the charging and discharging rate limits as follows:
\begin{align}\label{con:es_dynamics3}
& 0 \leq p_{k, t}^{\text{ch}} \leq \overline{P}_{k}^{\text{ch}} u_{k, t}^{\text{ch}}, 0 \leq p_{k, t}^{\text{dis}} \leq \overline{P}_{k}^{\text{dis}} u_{k, t}^{\text{dis}}, \forall k, t \nonumber \\
& 0 \leq u_{k, t}^{\text{ch}} + u_{k, t}^{\text{dis}} \leq 1, \forall k, t,
\end{align}
where $\overline{P}_{k}^{\text{ch}}$ and $\overline{P}_{k}^{\text{dis}}$ denote the maximum charged and discharged energy for the $k$-th ESS, respectively. Variables $u_{k, t}^{\text{ch}}$ and $u_{k, t}^{\text{dis}}$ are binary variables indicating the charging and discharging decisions of the $k$-th ESS, which are mutually exclusive during a specific time interval.

\subsubsection{Network Constraints}
Furthermore, we denote critical power loads as $d_{i, t}$ that must be satisfied for each CB $i$. Therefore, we have the following power balance equation:
\begin{align}\label{con:p_balance}
& \sum_{k} \left(p_{k, t}^{\text{dis}} - p_{k, t}^{\text{ch}} \right) + \sum_i w_{i, t} + p_{\text{base}, t} +  p_{\text{peak}, t} \\
& = \sum_i d_{i, t} + \sum_{l} p_{l, t} + \sum_{i} P_{i, t}^{\text{hvac}} + \sum_{v} p_{v, t}^{\text{ch}}, \forall t. \nonumber
\end{align}
where $w_{i, t}$ denotes output of the renewable in each CB $i$. Variables $p_{\text{base}, t}$ and $p_{\text{peak}, t}$ represent real-time power buy from a retail electricity market with the base price and the peak price, respectively.

Additionally, we use $\overline{P}$ to represent the capacity limit on the tie-line connected with the distribution system. Moreover, in the peak hours, a DSO may inform central controllers of CBs about reducing certain amount of power delivery through the PCC (in order to protect substations and transformers). Thus, $I_{t} \in [0.8, 1]$ is adopted to model uncertain DR signals from a DSO during peak hours. Then, we have the following constraints on the grid-connected tie-line:
\begin{equation}\label{con:pcc}
p_{\text{base}, t} + p_{\text{peak},t} = p_{\text{total},t}, 0 \leq p_{\text{total},t} \leq \overline{P} I_{t}, \forall t.
\end{equation}

\subsection{Uncertainty Set}

\subsubsection{RES Power Output}
As mentioned in previous sections, the robust optimization can handle various uncertainties through a pre-defined deterministic interval, such as $[\bar{w}_t - \hat{w}_t^-, \bar{w}_t + \hat{w}_t^+]$ for the RES power output $w$. Variable $\bar{w}_t$ represents the nominal value (i.e., forecasted value) of RES power output in a time slot $t$. Variables $\hat{w}_t^-$ and $\hat{w}_t^+$ are the maximum negative and positive deviations of the RES generation in a time slot $t$, respectively. Coordinating with the budget of uncertainty $\Gamma_w$, we have the cardinality uncertainty set for the RES generation as follows:
\vspace{-7pt}
\begin{align}
& \mathcal{W} := \left\{\mathbf{w} : w_{i, t} = \bar{w}_{i, t} + \hat{w}_{i, t}^+ u_{i, t}^+ - \hat{w}_{i, t}^- u_{i, t}^-, \forall i, t \right\}, \label{eq:w_set1} \\
& U := \bigg\{\boldsymbol{u}: \sum_{t = 1}^{T}\left(u_{i, t}^- + u_{i, t}^+ \right) \leq \Gamma_w, 0 \leq u_{i, t}^+, u_{i, t}^- \leq 1, \forall i, t \bigg\}. \nonumber  \vspace{-7pt}
\end{align}
where $u_{i, t}^+$ and $u_{i, t}^-$ are auxiliary variables indicating the degree of positive and negative deviation from the forecasted value $\bar{w}_{i, t}$.

\subsubsection{Load Demand}
Similarly, we can construct the uncertainty set for critical load demand as follows: 
\vspace{-7pt}
\begin{align}
& \mathcal{D} := \left\{\mathbf{d} : d_{i, t} = \bar{d}_{i, t} + \hat{d}_{i, t}^+ v_{i, t}^+ - \hat{d}_{i, t}^- v_{i, t}^-, \forall i, t \right\}, \label{eq:l_set1} \\
& V := \bigg\{\boldsymbol{v}: \sum_{t = 1}^{T}\left(v_{i, t}^- + v_{i, t}^+ \right) \leq \Gamma_d, 0 \leq v_{i, t}^+, v_{i, t}^- \leq 1, \forall i, t \bigg\}. \nonumber  \vspace{-7pt}
\end{align}
where $v_{i, t}^+$ and $v_{i, t}^-$ are auxiliary variables indicating the degree of positive and negative deviation from the forecasted value $\bar{d}_{i, t}$. Variable $\Gamma_d$ is the budget of uncertainty for elastic load demand.

\subsubsection{PEVs}
Moreover, the initial SoCs of PEVs are unknown, therefore, we have the following uncertainty set for PEV's SoC conditions:   
\vspace{-7pt}
\begin{align} \label{eq:z_set1}
& \mathcal{SoC} := \{\mathbf{SoC} : SoC_{v, 0} = \bar{SoC}_{v, 0} + \hat{SoC}_{v, 0}^+ z_{v, 0}^+  \\ \nonumber
& - \hat{SoC}_{v, 0}^- z_{v, 0}^-, \forall v, \boldsymbol{z} \in Z \}, \\
& Z := \bigg\{\boldsymbol{z}: \sum_{v = 1}^{N_v} \left(z_{v, 0}^- + z_{v, 0}^+ \right) \leq \Gamma_z, 0 \leq z_{v, 0}^+, z_{v, 0}^- \leq 1, \forall v \bigg\}. \nonumber \vspace{-7pt}
\end{align}
where $z_{v, 0}^+$ and $z_{v, 0}^-$ are auxiliary variables indicating the degree of positive and negative deviation from the forecasted value $\bar{SoC}_{v, 0}$. Variable $\Gamma_z$ is the budget of uncertainty for PEVs' initial SoCs.

\subsubsection{DR signals}
Additionally, we have the following uncertainty set for peak hours DR signals from the DSO:
\vspace{-7pt}
\begin{equation}\label{eq:I_set}
\mathcal{I} := \bigg\{\mathbf{I} : \sum_{t=1}^T (1-I_t) \leq \Gamma_I, 0.8 \leq I_t \leq 1, \forall t \bigg\}. \vspace{-5pt}
\end{equation}
where $\Gamma_I$ is the budget of uncertainty for a DSO's DR signals.

\subsection{Objective Function}
Our objective is to maximize consumers' comfort levels, while minimizing O\&M costs of the proposed campus-based commercial building system under the worst-case scenario. Thus, we formulated the energy management problem into a max-min-max problem as follows:
\begin{align}\label{con:primal}
\max & \left( \sum_t \sum_i J_{i, t} + \frac{N_v}{N_p} \sum_t \sum_v J_{v, t} + \sum_t \sum_l J_{l, t}\right) \kappa \nonumber \\
& + \min_{\boldsymbol{w}, \boldsymbol{d}, \boldsymbol{SoC}, \boldsymbol{I}} \max \bigg\{ - \sum_t \sum_k C_{\text{ESS}} ( p_{k, t}^{\text{ch}} + p_{k, t}^{\text{dis}} ) \nonumber \\
& \qquad \qquad - \sum_t \sum_v C_{\text{PEV}} p_{v, t}^{\text{ch}} \\
& \qquad \qquad - \sum_t C_{\text{base}} p_{\text{base}, t} - \sum_t C_{\text{peak}} p_{\text{peak}, t} \bigg\}.\nonumber
\end{align}
Note that $C_{\text{base}}$ and $C_{\text{peak}}$ are electricity prices for the base load and the peak load, respectively. The electricity price for the peak load is usually several times of the base load, which should be listed in the bilateral trading agreement that is signed by central controllers of CBs and a local utility company. Variables $C_{\text{ESS}}$ and $C_{\text{PEV}}$ are degradation cost coefficients for ESSs and PEVs, respectively. Variable $\kappa$ is the weighting factor representing the trade-off between O\&M costs and comfort levels.

\section{Solution Algorithm}\label{section3}

\subsection{Problem Reformulation}
The formulation above is a max-min-max problem, which cannot be solved directly by commercial solvers. Thus, as strong duality holds, we transform the inner $\max$ formulation into its dual as a $\min$ formulation in a matrix form, then using Big-M method to linearize nonlinear terms in the objective function of the new minimization problem:
\begin{align}\label{eq:dual1}
\min_{\boldsymbol{\lambda}, \boldsymbol{\mu}, \boldsymbol{\gamma}, \boldsymbol{\pi}, \boldsymbol{\sigma}} \quad &  - \boldsymbol{\lambda}^T \boldsymbol{F} + \boldsymbol{\mu}^T \boldsymbol{SoC} + \boldsymbol{\gamma}^T \boldsymbol{G} + \boldsymbol{\pi}^T \boldsymbol{w} + \boldsymbol{\sigma}^T \boldsymbol{d} \nonumber \\
\text{s.t.} \qquad & - \boldsymbol{\lambda}^T \boldsymbol{A} + \boldsymbol{\gamma}^T \boldsymbol{B} + \boldsymbol{H} \geq 0 \nonumber \\
\qquad & \boldsymbol{\mu}^T - M \boldsymbol{z} \leq 0, \boldsymbol{\mu} \pm \boldsymbol{\lambda} - M (1 - \boldsymbol{z} ) \leq 0 \nonumber \\
\qquad & \boldsymbol{\pi}^T - M \boldsymbol{u} \leq 0, \boldsymbol{\pi} \pm \boldsymbol{\lambda} - M (1 - \boldsymbol{u} ) \leq 0 \nonumber \\
\qquad & \boldsymbol{\sigma}^T - M \boldsymbol{v} \leq 0, \boldsymbol{\sigma} \pm \boldsymbol{\lambda} - M (1 - \boldsymbol{v} ) \leq 0 \nonumber \\
\qquad & \boldsymbol{\lambda}, \boldsymbol{\mu}, \boldsymbol{\gamma}, \boldsymbol{\pi}, \boldsymbol{\sigma} \in [0, 1],
\end{align}
subject to constraints~\eqref{eq:w_set1}--\eqref{eq:I_set}. Variables $\boldsymbol{\lambda}$, $\boldsymbol{\mu}$, and $\boldsymbol{\gamma}$ denote Lagrangian multipliers of constraints [\eqref{con:HVAC1}, \eqref{con:pev_dynamics1}, \eqref{con:EWH3}, \eqref{con:es_dynamics1}, \eqref{con:p_balance}, \eqref{con:pcc}], [\eqref{con:pev_dynamics2}, \eqref{con:es_dynamics2}--\eqref{con:es_dynamics3}], and [\eqref{con:HVAC2}, \eqref{con:EWH4}], respectively. Variables $\boldsymbol{F}$, $\boldsymbol{SoC}$, and $\boldsymbol{G}$ are sets of parameters in the constraints represented by $\boldsymbol{\lambda}$, $\boldsymbol{\mu}$, and $\boldsymbol{\gamma}$. Variables $\boldsymbol{A}$ and $\boldsymbol{B}$ are sets of parameters in the inner $\max$ problem~\eqref{con:primal}. However, there are two non-linear terms, $\boldsymbol{\pi}^T \boldsymbol{w}$ and $\boldsymbol{\sigma}^T \boldsymbol{d}$, in equation~\eqref{eq:dual1}, which makes the second-stage problem hard to solve. Thus, we need to linearize the non-linear terms through the Big-M method~\cite{LiCh17}.

\subsection{C\&CG Algorithm}
Combined with all the formulations from previous sections, we can finally reformulate the original problem into the following mixed-integer linear programming (MILP) problem.
\begin{align}\label{eq:ccg}
\max \quad &  \boldsymbol{\kappa}^T \boldsymbol{J} - \theta \nonumber \\
\text{s.t.} \qquad & \theta \geq - \boldsymbol{\lambda}^T \boldsymbol{F} + \boldsymbol{\mu}^T \boldsymbol{SoC} + \boldsymbol{\gamma}^T \boldsymbol{G} + \boldsymbol{\pi}^T \boldsymbol{w}+ \boldsymbol{\sigma}^T \boldsymbol{d} \nonumber \\
\qquad & constraints~\eqref{con:HVAC1}-\eqref{con:pcc}.
\end{align}
Variable $\boldsymbol{J}$ represents the constraints~\eqref{con:comfort1},~\eqref{con:comfort2}, and~\eqref{con:comfort3} that are related with comfort levels. As mentioned in previous sections, the worst-case scenario of each set of uncertainty is independent from all others and can only occur at a set's upper and lower bound, which are finite. To solve the problem in a reasonable solution time, we employ the C\&CG algorithm, while details are shown as follows:
\begin{enumerate}
    \item Set feasible first-stage decisions with sets of budget of uncertainties. Then, solve the sub problem~\eqref{eq:dual1} with the initial conditions. After that, use the answer of $\boldsymbol{u}_1$, $\boldsymbol{v}_1$, $\boldsymbol{z}_1$, and $\boldsymbol{I}_1$ as the worst-case scenario. Let $\text{UB} = +\infty$ be the upper bound, $\text{LB} = -\infty$ be the lower bound, and $1$ be the iteration number.

    \item Solve the master problem~\eqref{eq:ccg} with the worst-case scenario. Determine the updated optimal first-stage decision and budget of uncertainty. Let $\text{LB} =  \kappa \boldsymbol{J}_s - \theta$.

    \item Solve the sub problem~\eqref{eq:dual1} with the updated optimal solutions from Step 2. Make $\text{UB} = \min\{\text{UB}, \kappa \boldsymbol{J}_s - \rho^*\}$, where $\rho^*$ is the updated value of the objective function.

    \item If $\text{UB} - \text{LB} \le \varepsilon $, then stop. Otherwise, update $s$ to $s + 1$ and return to Step 2.
\end{enumerate}

\section{Simulation Results}\label{section4}
We set $10^{-2}$ as the convergence tolerance. All simulations are implemented on a desktop with a 3.00 GHz Intel Core i5-7400 CPU and 8GB RAM.

\subsection{Numerical Settings}
\begin{table}[!ht]
\centering
\vspace{-15pt}
\tabcolsep=0.11cm
\caption{Comfort Level Related Parameters}
\label{table:comfortlevel}
\begin{tabular}{cccc}
\cline{1-4}
Type & $T^d_{i}$ ($^\circ C$) & $\delta_{i}$ ($^\circ C$) & $\epsilon_{i}$ ($^\circ C$) \\
\hline
HVAC & 24 & 2 & 0.5  \\
\hline
Type & $SoC_{v}^d$ ($\%$) & $SoC_{v}^{\text{base}}$ ($\%$) & $SoC_{v, 0}$ ($\%$) \\
\hline
PEV & 80 & 10 & Uncertain \\
\hline
Type & $T_{l}^d$ ($^\circ C$) & $\delta_{l}$ ($^\circ C$) & $T_{l, 0}$ ($^\circ C$)\\
\hline
EWH & 40 & 10 & 30 \\
\hline
\end{tabular}
\end{table}
\begin{table}[!ht]
\centering
\vspace{-15pt}
\tabcolsep=0.11cm
\caption{ESS/PEV Parameters}
\label{table:3}
\begin{tabular}{cccc}
\cline{1-4}
$SoC_{k, 0}$ ($\%$) & $\underline{SoC}_k$/$\overline{SoC}_k$ ($\%$) & $\overline{E}_{k}$ (kWh) & $\overline{P}_{k}^{\text{ch}}$/$\overline{P}_k^{\text{dis}}$ (kW) \\
\hline
50 & 5/95 & 80 & 4/4    \\
\hline
Type & $\underline{SoC}_v$/$\overline{SoC}_v$ ($\%$) & $\overline{E}_{v}$ (kWh) & $\overline{P}_{v}^{\text{ch}}$ (kW) \\
\hline
Tesla Model S & & & \\
75D & 5/95 & 75 & 11.5 \\
\hline
Tesla Model X & & & \\
100D & 5/95 & 100 & 17.2 \\
\hline
Nissan Leaf & & & \\
SV & 5/95 & 30 & 3.6 \\
\hline
\end{tabular}
\end{table}
The considered commercial campus consists of six CBs, six packs of roof-top solar panels, fifty PEVs in three different types, and one aggregated critical power load. Each CB includes one EWH, one HVAC system, and one pack of batteries, which are scalable. The comfort levels are only considered during the office hours when the consumers are on the commercial campus (from $8$ a.m. to $8$ p.m.), where the parameters are from a real-world commercial building, as shown in Tables~\ref{table:comfortlevel},~\ref{table:3}. The total installed capacity of the solar panels is $60 kW$, where the historic generation patterns are taken from~\cite{Wind}. Also, the historical data of the solar irradiance and the outdoor temperature are from~\cite{Wind}, with proper scaling coefficients. The base and peak electricity prices are taken from a real-world utility company, where the peak electricity price is 10 times of the base electricity price. The EWHs' power to heat ratio is set to $1.2$. The data for EWHs are taken from~\cite{SaOr15}. The electricity exchange limit between the proposed system and the retail electricity market is set to be $1867 kW$. The degradation cost coefficients $C_\text{ESS}$ and $C_{\text{PEV}}$ are set to be $0.0035 \$/kWh$.

\subsection{Simulation Results}
\begin{figure}[!t]
\centering
\begin{minipage}[b]{0.45\linewidth}
\includegraphics[width = 1.65in]{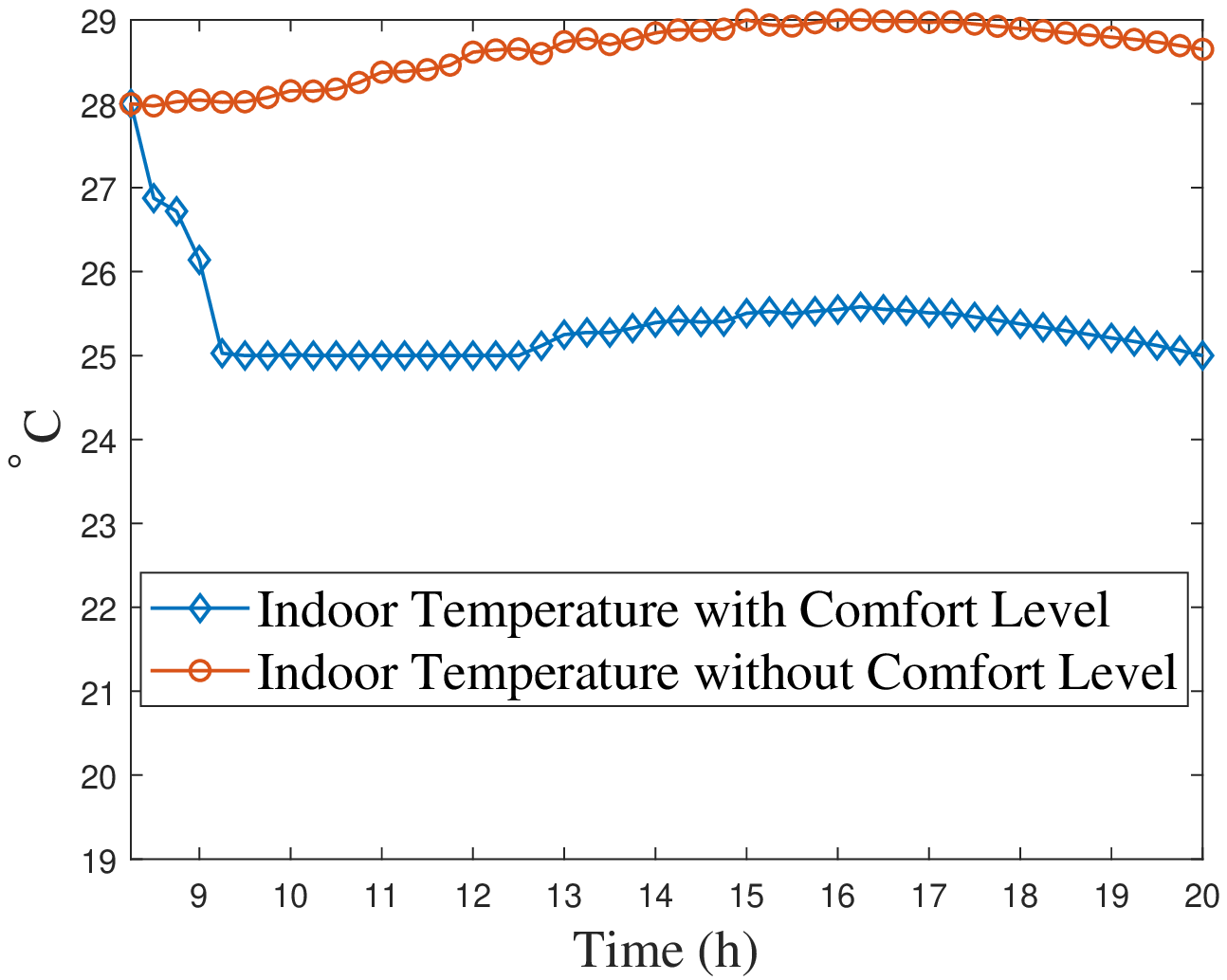}
\caption{The indoor temperature with/without the proposed comfort level model.}
\label{fig:temp}
\end{minipage}
\quad
\begin{minipage}[b]{0.45\linewidth}
\includegraphics[width = 1.65in]{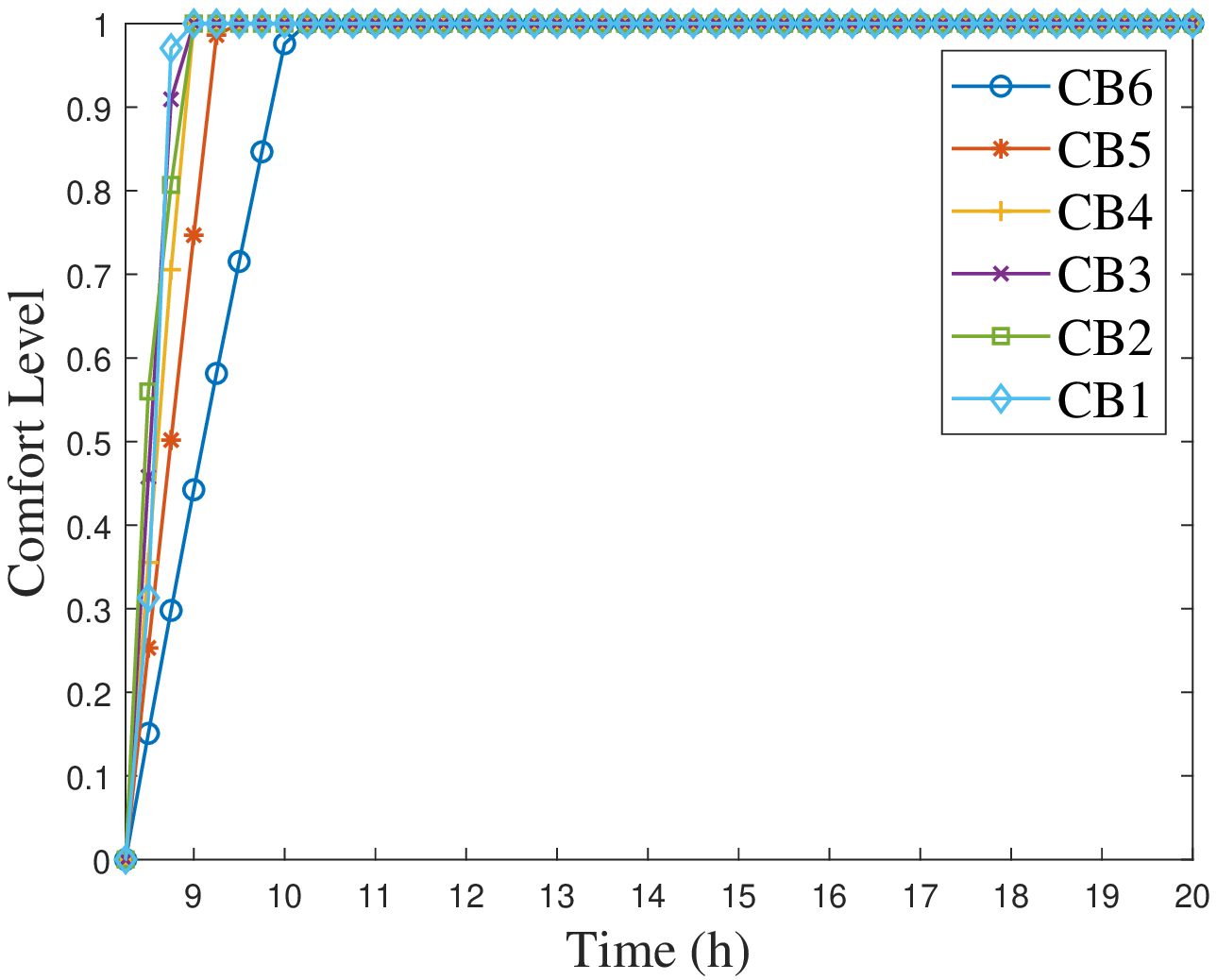}
\caption{The comfort level related HVAC systems considering the proposed comfort level model.}
\label{fig:HVAC_D}
\end{minipage}
\end{figure}
\begin{figure}[!t]
\centering
\vspace{-15pt}
\begin{minipage}[b]{0.45\linewidth}
\includegraphics[width = 1.65in]{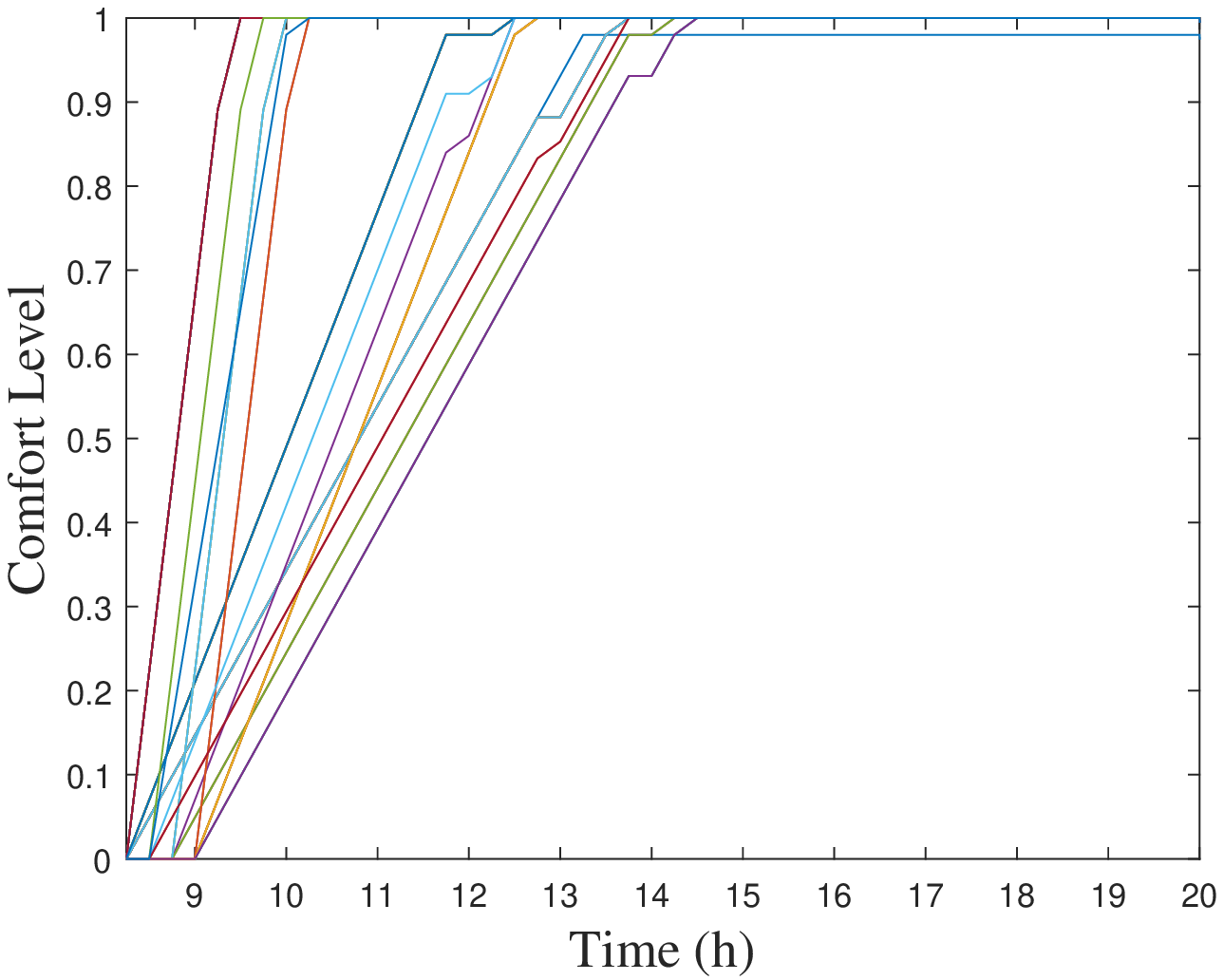}
\caption{The comfort level related with PEVs considering the proposed comfort level model.}
\label{fig:PEV_D}
\end{minipage}
\quad
\begin{minipage}[b]{0.45\linewidth}
\includegraphics[width = 1.65in]{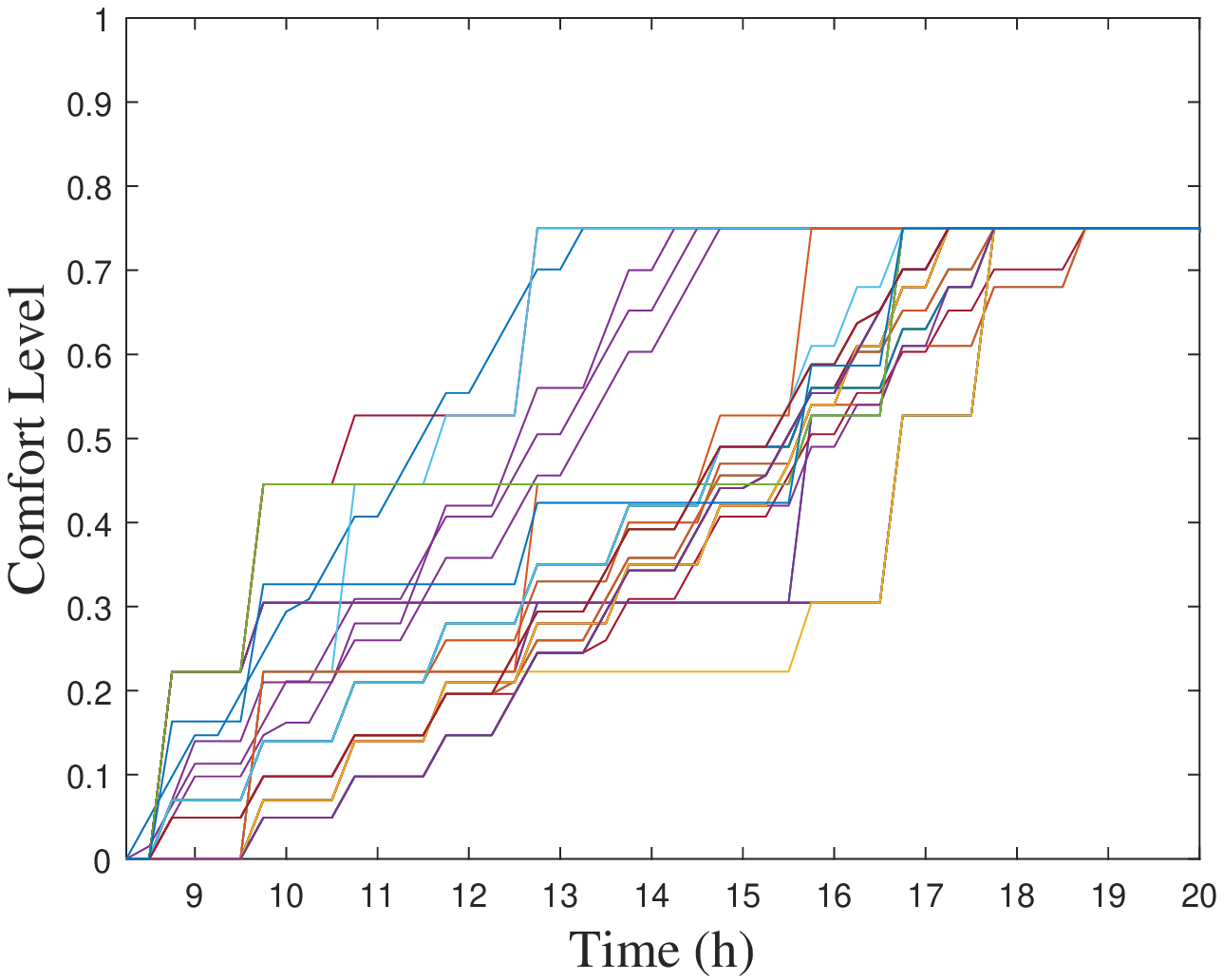}
\caption{The comfort level related with PEVs without the proposed comfort level model.}
\label{fig:PEV_DnoJ}
\end{minipage}
\end{figure}
\begin{figure}[!t]
\centering
\vspace{-15pt}
\begin{minipage}[b]{0.45\linewidth}
\includegraphics[width = 1.65in]{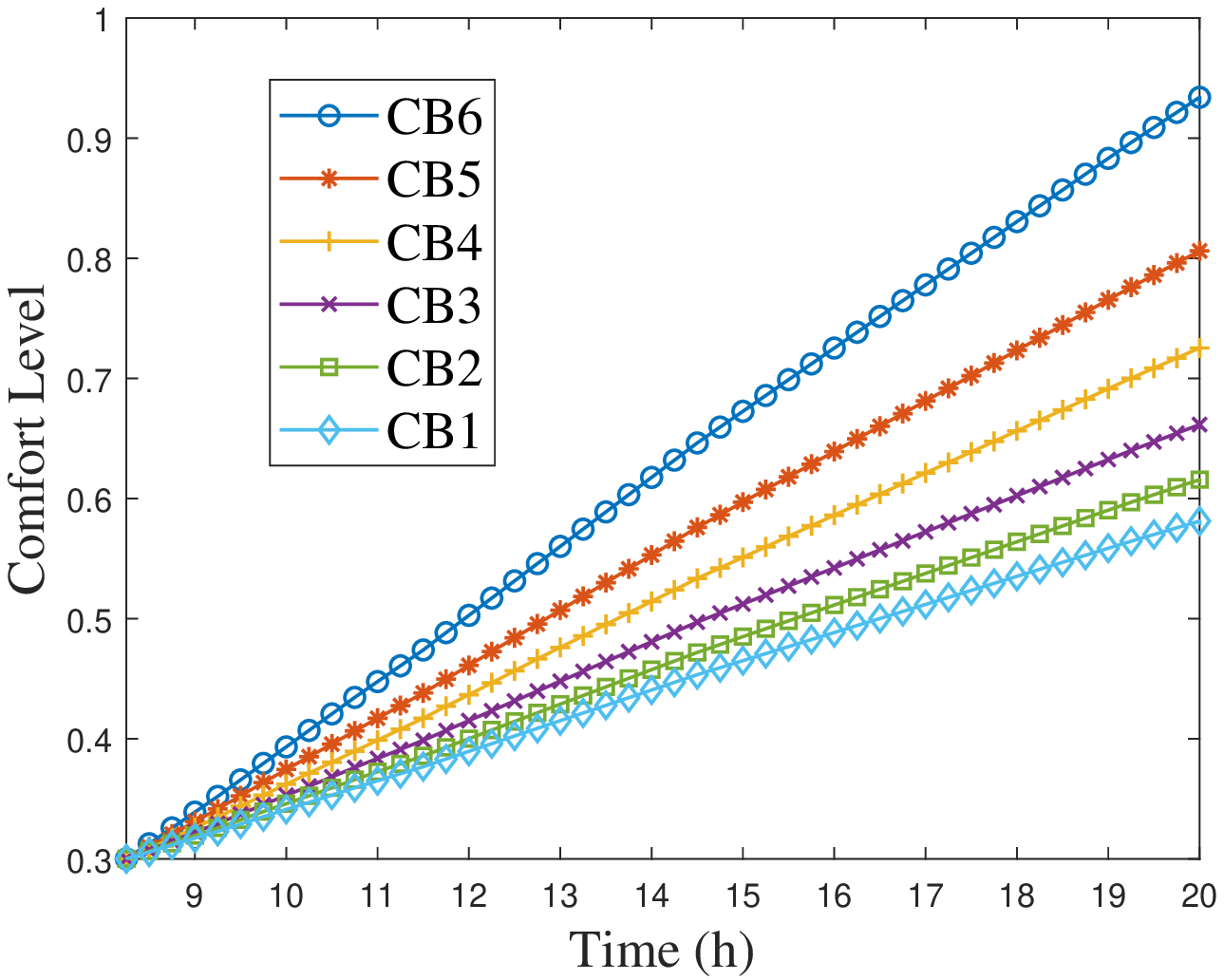}
\caption{The comfort level related with EWHs considering the proposed comfort level model.}
\label{fig:EWH_D}
\end{minipage}
\quad
\begin{minipage}[b]{0.45\linewidth}
\includegraphics[width = 1.65in]{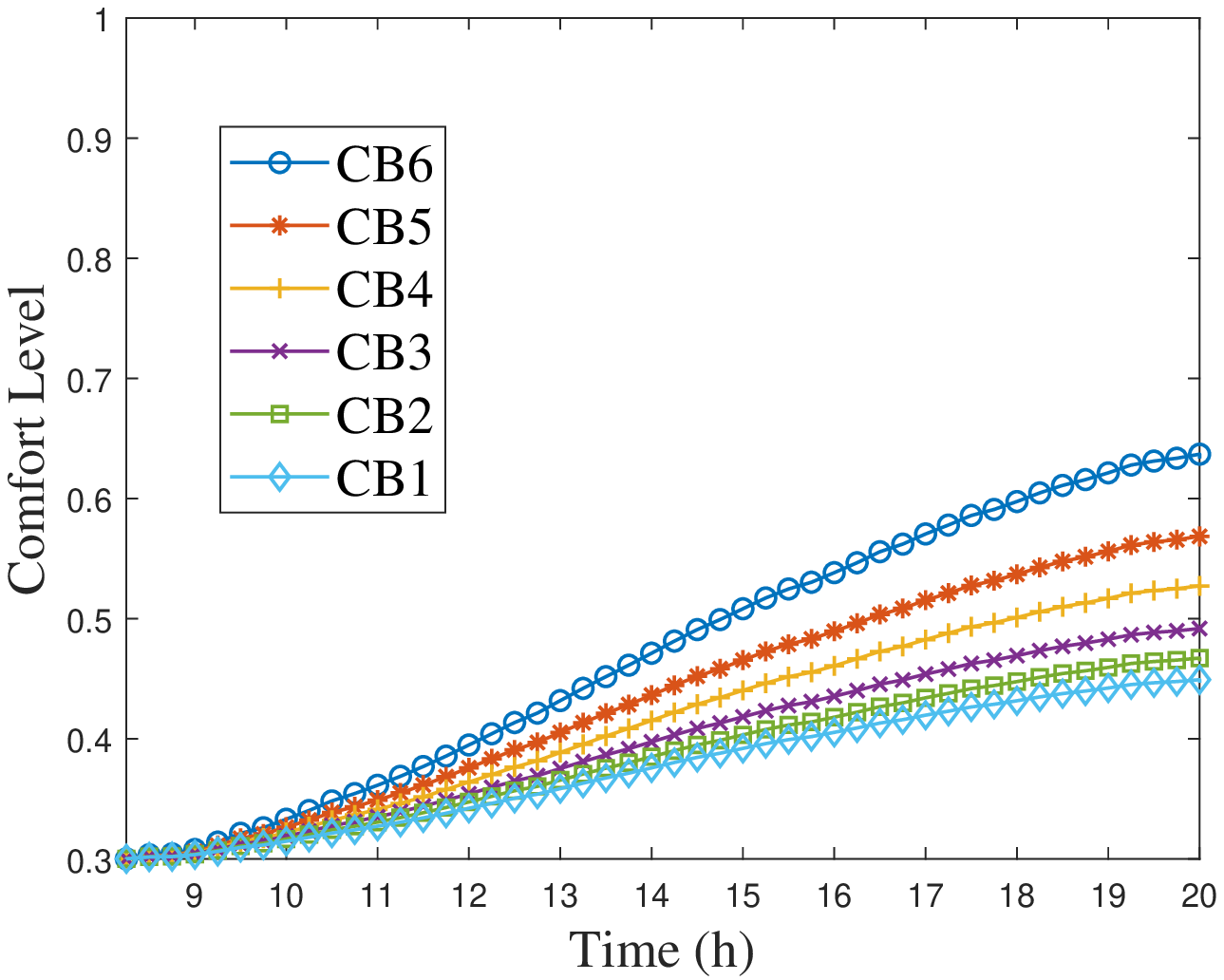}
\caption{The comfort level related with EWHs without the proposed comfort level model.}
\label{fig:EWH_DnoJ}
\end{minipage}
\end{figure}
In order to show the effectiveness of the proposed comfort level models, we test the system in a single scenario with/without comfort level related objectives, while the comfort levels are still calculated through the original constraints in the system modeling section. As shown in Fig.~\ref{fig:temp}, the indoor temperatures of the first building reduces significantly when considering the proposed comfort level formulations, which is because the base electricity price during this period (from 8 am to 10 am) is lower enough to provide sufficient power to the HVAC system to reduce the indoor temperature. This situation is also demonstrated in Fig.~\ref{fig:HVAC_D}, where the comfort level increases vastly during the same period.

We further test our comfort level models of PEVs and EWHs in the same circumstances with/without comfort levels related objective functions. As shown in Figs.~\ref{fig:PEV_D} and~\ref{fig:PEV_DnoJ}, comfort levels increase faster when considering the proposed comfort level model. The PEVs are also charged to the desired energy level when the comfort level model is not in use. However, in this situation, only the electricity price impact is taken into consideration, which may result in less comfort when the consumer left earlier.
\balance
Moreover, we test the influence of comfort level on EWHs. As shown in Figs.~\ref{fig:EWH_D} and~\ref{fig:EWH_DnoJ}, the comfort level for EWHs with the proposed model are much better. This is because EWHs can provide most of the hot water demand as needed. Besides, through checking the trade-off between electricity price and comfort level, the comfort level of the case with the proposed model can still maintain the comfort level of EWHs at a high level compared with the one without the proposed model.

The overall comfort levels increased from $0.42$ to $1$ while the total cost only increases $19.9\%$. Therefore, our proposed comfort level model for HVAC systems, PEVs and EWHs are more reasonable and suitable for the commercial building system.

\section{Conclusion}\label{section5}
In this paper, we propose several novel comfort level models to minimize O\&M costs of campus-based CBs while maximize various comfort levels simultaneously under worst-case scenarios. ARO has been leveraged to handle various uncertainties within the proposed system: (i) demand response signals from the DSO, (ii) arrival SoC conditions of PEVs, (iii) power outputs of RES units, and (iv) load demand of other appliances. A C\&CG algorithm is proposed to solve the reformulated NP-hard min-max problem. Extensive simulation results have demonstrated the effectiveness of our proposed optimal energy management strategy for campus-based CBs in both minimizing O\&M costs and maximizing comprehensive comfort levels.

\bibliographystyle{IEEEtran}
\bibliography{./mybibfile}

\end{document}